\newtheorem{lemma}{Lemma}[section]
\newtheorem{prop}[lemma]{Proposition}
\newtheorem{theorem}[lemma]{Theorem}

\newtheorem{rem}[lemma]{Remark}

\newtheorem{frem}[lemma]{Final remark}

\newcommand{\kla}{\left ( }
\newcommand{\nach}{\rightarrow}
\newcommand{\mer}{\right ) }

\newcommand{\for}{\begin{eqnarray*}}
\newcommand{\mel}{\end{eqnarray*}}
\newcommand{\qed}{\hspace*{\fill}$\Box$\hz\pagebreak[1]}

\newcommand{\mitt}{\left | { \atop } \right.}
\newcommand{\kl}{\pl \le \pl}
\newcommand{\gl}{\pl \ge \pl}
\newcommand{\lel}{\pl = \pl}
\newcommand{\lell}{\p=\p}
\newcommand{\kll}{\p \le \p}
\newcommand{\gll}{\p \ge \p}

\newcommand{\nz}{{\rm  I\! N}}
\newcommand{\nen}{n \in \nz}

\newcommand{\ez}{{\rm I\! E}}
\newcommand{\rz}{{\rm  I\! R}}

\newcommand{\p}{\hspace{.05cm}}
\newcommand{\pl}{\hspace{.1cm}}
\newcommand{\pll}{\hspace{.3cm}}
\newcommand{\pla}{\hspace{1.5cm}}
\newcommand{\hz}{\vspace{0.5cm}}

\newcommand{\si}{\sigma}

\newcommand{\eps}{\varepsilon}

\newcommand{\lzn}{\ell_2^n}

\newcommand{\lin}{\ell_{\infty}^n}
\newcommand{\lif}{\ell_{\infty}}

\newcommand{\lgi}{\ell_{g,\infty}}

\newcommand{\lgx}{\ell_{g,\infty}^{(x)}}

\newcommand{\Lgx}{ {\cal L}_{g,\infty}^{(x)}}

\renewcommand{\L}{{\cal L}}

\newcommand{\ef}{{\cal F}}

\newcommand{\noo}{\left \|}
\newcommand{\tno}{\big\|\!\big|}
\newcommand{\rrm}{\right \|}

\newcommand{\bet}{\left |}
\newcommand{\rag}{\right |}

\newcommand{\summ}{\sum\limits}

\newcommand{\su}{\subset}

\newcommand{\pq}{\pi_{pq}}

\documentstyle[11pt]{article}
\begin{document}

\begin{center}

{\bf \huge Cotype and summing properties in Banach spaces}\hz

{\large Marius Junge}
\end{center}
\hz\hz
\newtheorem{llemma}{Lemma}
\newtheorem{pprop}[llemma]{Proposition}
\newtheorem{heorem}[llemma]{Theorem}

\begin{abstract}
\parindent0em
\hspace{-1,9em} It is well known in Banach space theory that for a
finite
dimensional space $E$ there exists a constant $c_E$, such that for all
sequences $(x_k)_k \subset E$ one has
\[ \summ_k \noo x_k \rrm \kl c_E \pl \sup_{\eps_k \pm 1} \noo \summ_k
\eps_k x_k \rrm \pl .\]
Moreover, if $E$ is of dimension $n$ the constant $c_E$ ranges between
$\sqrt{n}$ and $n$. This implies that absolute convergence and
unconditional
convergence only coincide in finite dimensional spaces. We will
characterize Banach spaces $X$, where the constant $c_E \sim \sqrt{n}$
for all finite dimensional subspaces.
More generally, we prove that an estimate $c_E \kll c
n^{1-\frac{1}{q}}$
holds for all $n \in \nz$ and all $n$-dimensional subspaces $E$ of $X$
if and only
if the eigenvalues of every operator factoring through
$\ell_{\infty}$ decrease of order $k^{-\frac{1}{q}}$ if and only if
$X$ is of weak cotype $q$, introduced by Pisier and Mascioni.
We emphasize that in contrast to Talagrand's equivalence
theorem on cotype $q$ and absolutely $(q,1)$-summing spaces this
extends
to the case $q=2$. If $q>2$  and one of the conditions above is
satisfied
one has
\[ \kla \summ_k \noo x_k \rrm^q \mer^{\frac{1}{q}} \kl C^{1+l}\pl
(1+{\rm log}_2)^{(l)}((1 +{\rm log}_2 n)^{\frac{1}{q}}) \pl \ez \noo
\summ_k \eps_k x_k \rrm \]
for all $n,l \in \nz$ and $(x_k)_k \subset E$, $E$ a $n$ dimensional
subspace of $X$.
In the case $q=2$ the same holds if we replace the expected value by
the supremum.
\end{abstract}

\section*{Introduction}
\setcounter{lemma}{0}

In  Banach spaces unconditional
convergence and absolute convergence only coincide for finite
dimensional
spaces. More precisely, a constant $0<c<\infty$ such that
\[ \summ_k \noo x_k \rrm \kl c \pl \sup_{\eps_k\pm 1} \noo \summ_k
\eps_k x_k \rrm \]
holds for all sequences $(x_k)_k\subset X$ if and only if $X$ is of
finite dimension.
The best possible constant $c$ is called the absolutely $1$-summing
norm of
the identity of $X$ ($\pi_1(Id_X)$). This notion was originally
introduced by
Grothendieck under the name 'semi-integrale \'{a} droite'. But Orlicz
discovered before that unconditional converging series are al least
absolutely
$2$-summing, provided the underlaying spaces is $L_p$, $1\le p \le 2$.
That's why this property is called Orlicz property. It is best
possible,
since Dvoretzky's theorem
ensures that for each $\delta>0$, $\nen$ there are elements
$(x_k)_{k=1}^n$ in an infinite dimensional Banach space $X$ which
satisfies
\[ \summ_1^n \noo x_k \rrm \gl (1-\delta) \pl \sqrt{n} \quad \mbox{and}
\quad
\sup_{\bet \alpha_k\rag \le 1} \noo \summ_1^n \alpha_k x_k \rrm \kl 1
\pl .\]
We will study spaces where this estimate is optimal or a certain
growth rate occurs. This is contained in the following
\pagebreak

\begin{heorem} Let $2\kll q <\infty$. For a complex Banach space X the
following properties are equivalent.
\begin{enumerate}
\item[i)] There exists a constant $c>0$ such that for all $\nen$ and
all $n$ dimensional subspace
$E\subset X$ one has
\[ \summ_k\noo x_k \rrm \kl c\pl n^{1-\frac{1}{q}} \pl
\sup_{\eps_k=\pm 1}\noo \summ_1^n \eps_k x_k \rrm\]
for all sequences $(x_k)_k \subset E$.
\item[ii)] $X$ is of weak cotype $q$, in other words there exists a
constant $0< c_2 < \frac{1}{8e}$ such that for all $\nen$ and
$x_1,..,x_n$
\[\left .
\begin{array}{c}
\summ_1^n \bet \langle x_k,x^* \rangle \rag^2 \kll \noo x^* \rrm\\
\mbox{and} \quad \noo x_k \rrm \gl c_2\quad \mbox{for} \quad k=1,..,n\\
\end{array}\right \}
\quad \pl \Rightarrow \quad\pl \ez \noo \summ_1^n \eps_k x_k \rrm  \gl
c_2 \pl n^{\frac{1}{q}}\pl .\]
\item[iii)] There exists a constant $c_3$ such that for all operators
$T :X \nach X$ which factors through $\ell_{\infty}$, i.e.
$T\p=\p SR$, $R:X \nach \ell_{\infty}$ and $S: \ell_{\infty} \nach X$
one has
\[ \sup_{\nen} n^{\frac{1}{q}}|\lambda_n(T)| \kl c_3 \pl \noo S \rrm
\pl \noo R \rrm \pl,\]
where $(\lambda_n(T))_{\nen}$ denotes the eigenvalue sequence of $T$ in
non-increasing order
according to there multiplicity.
\end{enumerate}
If $q>2$ and one of the conditions above are satisfied there is a
constant $C$ such that
\[ \kla \summ_k\noo x_k \rrm^q \mer^{\frac{1}{q}} \kl C^{1+l} \pl
(\max\{1,{\rm log}_2\})^{(l)}((1 +{\rm log}_2 n)^{\frac{1}{q}}) \pl \ez
\noo \summ_k \eps_k x_k \rrm \]
holds for all $n,l \in \nz$ and $(x_k)_k \subset E$, $E$ a
$n$-dimensional
subspace of $X$.
\end{heorem}

This theorem is somehow at end of a fruitful investigation of summing
and cotype
properties in Banach spaces. Starting point is certainly the pioneering
work of Maurey and Pisier \cite{MP}. In their paper they obtained the
equivalence
in terms of the cotype index. Using deep methods from the theory of
stochastic
processes, the so-called concentration phenomena, Talagrand improved
Maurey/Pisier's result.

\begin{heorem}[Talagrand] Let $2<q<\infty$ and $X$ a Banach space the
following are equivalent
\begin{enumerate}
\item The identity of $X$ is absolutely $(q,1)$-summing, i.e. there is
a constant $c_1>0$ such that
\[ \kla \summ_k \noo x_k \rrm^q \mer^{\frac{1}{q}} \kl c_1 \pl
\sup_{\eps_k \pm 1} \noo \summ_k \eps_k x_k \rrm \pl,\]
for all $(x_k)_k\subset X$.
\item $X$ is of cotype $q$, i.e. there is a constant $c_2$ such that
\[ \kla \summ_k \noo x_k \rrm^q \mer^{\frac{1}{q}} \kl c_2 \pl  \ez
\noo \summ_k \eps_k x_k \rrm \pl,\]
for all $(x_k)_k\subset  X$.
\end{enumerate}
\end{heorem}

In Talagrand's theorem the case $q=2$ is not included and for a good
reason:
\pagebreak

\begin{heorem}[Talagrand] There is a symmetric sequence space which has
the Orlicz property but is not of cotype $2$.
\end{heorem}

Nevertheless, in the proof of the main theorem,
also in the case $q=2$, we heavily use the probabilistic machinery
established
by Talagrand. For $q>2$ the modified cotype condition in $ii)$ can be
replaced by the usual cotype condition restricted to vectors of equal
norm.
This is not possible for $q=2$, since equal norm cotype $2$ is the same
as
cotype $2$. But this modified condition turns out to be a basic tool
for the application
of the probabilistic method.

By the way, using the main theorem Talagrand's example yields
a symmetric sequence space which is of weak cotype $2$ but not of
cotype $2$.
This is impossible in the category of weak Hilbert spaces, since every
symmetric
weak Hilbert space is actually a Hilbert space.

In this setting the 'weak' theory is more adapted
to prove abstract characterization theorems than the 'strong' theory.
This is also
true for eigenvalue
estimates. It happens quite often that eigenvalue estimates for weak
$\ell_p$
spaces are easier to prove than eigenvalue estimates for the spaces
$\ell_p$ themselves.
A useful tool in this context is notion of Weyl numbers. The connection
between Weyl numbers and weak cotype was actually discovered by
Mascioni \cite{MAS}.
We should note that the equivalence between eigenvalue estimates and
summing properties
can be proved using a generalization of Maurey's theorem, provided
$q>2$.
This approach was pursued in \cite{J1,J2}.

Finally we come to the estimate with the iterated logarithm. This will
be investigated in chapter 2 and is based on the introduction
of optimal cotype spaces. The idea is to measure cotype and summing
conditions
in terms of maximal, symmetric sequence spaces. It turns out
that a certain self concavity is a generalization of the
submultiplicativity conditions which occurred in the basic paper of
Maurey and
Pisier. This broader framework turns out to be more natural
to describe cotype conditions of Orlicz spaces,
although we will not start this investigation here. In order to
find the best possible eigenvalue behavior of operators factoring
through $\ell_{\infty}$ we will also proof the main theorem in
a slightly more general setting.

\setcounter{section}{0}
\section*{Preliminaries}
In what follows $c_0, c_1,$ .. always denote universal constants.
We use standard Banach space notation. In particular, the classical
spaces $\ell_q$ and $\ell_q^n$, $1\le q\le \infty$, $\nen$, are defined
in the usual way. We will also use the Lorentz spaces $\ell_{pq}$ where
$1\kll p,q \kll \infty$. This space consists of all sequences
$\si \in \lif$ such that

\[\noo \si \rrm_{pq} \pl
  := \pl \kla \summ_n \kla n^{\frac{1}{p}} \, \si_n^* \mer^q n^{-1}
  \mer^{\frac{1}{q}} \,<\,
  \infty \, \pll.   \]

For $q=\infty$ the needed modification is given by

\[ \noo \si \rrm_{p\infty}  \pl := \pl  \sup_{\nen} \,
n^{\frac{1}{p}}\,\si_n^*
   \pl <\pl \infty . \]

Here $\si^*\,=\,(\si_n^*)_{\nen}$ denotes the
non-increasing rearrangement of $\si$. More generally, for a non
decreasing
sequence $(g(n))_{\nen}$ with $g(1)=1$ we denote by $\ell_{g,\infty}$
the space of sequences $\si$ such that
\[ \noo \si \rrm_{g,\infty} \pl :=\pl  \sup_n g(n) \pl \si_n^* \pl <
\pl \infty \pl .\]

The standard reference on operator ideals is the monograph of Pietsch
\cite{PIE}. The ideals of linear bounded operators, finite rank
operators,
integral operators are denoted by $\L$, $\ef$.

Let $1 \le q \le p \le \infty$ and $\nen$. For an operator $T \in
\L(X,Y)$
the pq-summing norm of T with respect to $n$ vectors is defined by
\[ \pq^n(T) \pl := \pl
 \sup\left\{\p \kla \summ_1^n \noo Tx_k \rrm^p \mer^{1/p} \p \bet \pl
  \sup_{\noo x* \rrm_{X^*}\le1} \kla \summ_1^n \bet\langle
  x_k,x^*\rangle \rag^q \mer^{1/q}
 \pl \le \pl 1\right.\p \right\} \pl .\]
An operator is said to be absolutely pq-summing, short pq-summing,
$(T \in \Pi_{pq}(X,Y))$ if
\[ \pq(T) \pl := \pl \sup_n \pq^n(T) \pl < \pl \infty \pl . \]
Then $(\Pi_{pq},\pq)$ is a maximal and injective Banach ideal (in the
sense
of Pietsch). As usual we abbreviate  $(\Pi_q,\pi_q) :=
(\Pi_{qq},\pi_{qq})$.
For further information about absolutely pq-summing operators we refer
to
the monograph of Tomczak-Jaegermann \cite{TOJ}.

In the following $(\eps_k)_{k\in \nz}$, $(g_k)_{k\in \nz}$ denotes a
sequence
of independent normalized Bernouilli, gaussian variables. A Banach
space
$X$ is of Rademacher, gaussian cotype $q$ if there exists a constant
$c>0$ such that
for all sequences $(x_k)_1^n \subset X$ one has
\for
\kla \summ_k \noo x_k \rrm^q \mer^{ \frac{1}{q} } &\le& c \pl \ez \noo
\summ_1^n \eps_k x_k \rrm\pl \mbox{or}\pl
\kla \summ_k \noo x_k \rrm^q \mer^{\frac{1}{q}} \kl c \pl \ez \noo
\summ_1^n g_k x_k \rrm\pl \mbox{resp.}
\mel
Here and in the following $\ez$ means expected value. The best possible
constant will by denoted by $Rc_q(X)\p:=\p Rc_q(id_X)$,
$c_q(X) \p:=\p c_q(Id_X)$, respectively. If this definition is
restricted to $n$ vectors
we write $RC_q^n$, $c_q^n$, respectively. As usual we will use the
abbreviation
\[ \ell(u) \pl:=\pl \sup_n \kla \ez \noo \summ_1^n g_k u(e_k) \rrm^2
\mer^{\frac{1}{2}} \]
for all operator $u \in \L(\ell_2,X)$. Here $(e_k)_k$ is the sequence
of unit vectors.
By the rotation invariance this norm is invariant by orthogonal
transformation
of this basis.

Finally some s-numbers are needed. For an operator $T\in\L (E,F)$ and
$\nen$
the $n$-th $approximation$ $number$ is defined by

\[ a_n(T) \pl :=\pl \inf\{\, \noo T-S \rrm \, | \,rank(S)\,< \, n \,\}
\pla ,\]

whereas the $n$-th $Weyl\,number$ is given by

\[ x_n(T) \pl :=\pl \sup\{\, a_n(Tu)\, |\, u \in \L(\ell_2,E) \,
\mbox{with}
 \, \noo u \rrm \,\le \, 1\,\} \pla .\]
Let $s \in \{a,x\}$. By $\L_{pq}^{(s)}$, $\L_{g,\infty}^{(S)}$ we
denote the ideal of operators $T$
such that $(s_n(T))_{\nen}\in \ell_{pq}$, $(s_n(T))_{\nen}\in
\ell_{g,\infty}$ with the associated quasi-norms
\[ \ell_{pq}^{(s)}(T) \pl:=\pl \noo (s_n(T))_{\nen} \rrm_{\ell_{pq}}
\quad
\mbox{and}\quad \ell_{g,\infty}^{(s)}(T) \pl:=\pl \noo (s_n(T))_{\nen}
\rrm_{\ell_{g,\infty}} \pl .\]

\section{Proof of the main theorem}
\setcounter{lemma}{0}

We will proof our main theorem in a little bit broader framework. For
the eigenvalue
estimate we allow a certain growth rate $(g(n))_{n \in \nz}$. Certainly
some reasonable conditions are required.

\begin{enumerate}
\item[S)]
\begin{enumerate}
\item[i)] $g(1)\lell1$ and $(g(n))_{n \in \nz}$ non decreasing.
\item[ii)] There exists a constant $S_2$ such that for all $1\kll k
\kll n$
\[ g(n) \kl S_2 \pl \frac{n}{k} \pl g(k) \pl .\]
\item[iii)] The space $\lgi$ is equivalent to a normed space. The
equivalence constant
is denoted by $S_3$.
\item[iv)] There is a constant $S_4$ such that for all $n \in \nz$
\[ \summ_1^n \frac{1}{g(k)} \kl S_4 \pl \frac{n}{g(n)} \pl .\]
\end{enumerate}
\item[L)] There exists a $0< t < \infty$ such that $n^{\frac{1}{t}} \kl
L_t \pl g(n)$.
\item[M)] There exists a natural number $r\gll 2$ with $t \kll r$ and a
constant $M_r$
such that
\[ g(k^{2r}) \kl M_r \pl g(k^r) \pl g(k)^r \pl .\]
\end{enumerate}

Condition $iii)$ and $iv)$ are actually equivalent but there is no need
to
go into further details. The last condition $M)$ is clearly satisfied
for supermultiplicative sequences.

\begin{theorem} Let $g$ be a sequence which satisfies the conditions
$Si)-iv)$, $L)$, and $M)$ and
set $D\p:=\p 2^{\frac{9}{2}}e^{\frac{3}{2}}S_2^2$.
For a complex Banach space X the following properties are equivalent.
\begin{enumerate}
\item[i)] There exists a constant $c_1>0$ such that for all
$n$-dimensional subspaces $E\subset X$ of $X$ one has
\[ \pi_1(Id_X) \kl c_1 \pl \frac{n}{g(n)} \pl .\]
\item[i')] There exists a constant $c'_1>0$ such that for all $\nen$
and $x_1,..,x_n \in X$ one has
\[ \summ_1^n\noo x_i \rrm \kl c'_1\pl \frac{n}{g(n)} \pl
\sup_{\eps_i=\pm 1}\noo \summ_1^n \eps_i x_i \rrm\pl .\]
\item[ii)] $X$ is of weak cotype $G$ or there exists a constant $c_2>0$
such that for all $\nen$ and $x_1,..,x_n$
\[\left .
\begin{array}{c}
\summ_1^n \bet \langle x_k,x^* \rangle \rag^2 \kll \noo x^*\rrm\\
\mbox{and} \quad \noo x_k \rrm \gl \frac{1}{D}\quad \mbox{for} \quad
k=1,..,n\\
\end{array}\right \}
\quad \pl \Rightarrow \quad\pl \ez \noo \summ_1^n \eps_k x_k \rrm \gl
\frac{1}{c_2} \pl g(n)\pl .\]
\item[iii)] There exists a constant $c_3$ such that for all operators
$T :X \nach X$ which factors through $\ell_{\infty}$, i.e.
$T\p=\p SR$, $R:X \nach \ell_{\infty}$ and $S: \ell_{\infty} \nach X$
one has
\[ \sup_{\nen}g(n)\pl |\lambda_n(T)| \kl c_3 \pl \noo S \rrm \pl \noo R
\rrm \pl,\]
where $(\lambda_n(T))_{\nen}$ denotes the eigenvalue sequence of $T$ in
non-increasing order
according to there multiplicity.
\end{enumerate}
\end{theorem}

For the proof of the main result we will closely follow Talagrand's
work.
The main difference occurs when we establish a situation where the
concentration phenomena can be applied. Let us recall
that an operator $T \in \L(X,Y)$ is of weak cotype $g$ ($T \in
WC_g(X,Y)$),
if there is a constant $c>0$ such that
\[ \sup_k g(k) \pl a_k(Tu) \kl c \pl \ell(u) \pl .\]
The norm $wc_g(T)$ is defined as the infimum over all $c$ satisfying
the inequality
above. The following lemma is well known and at origin of the so called
weak
theory, see \cite{PSW}. Nevertheless, we give a proof in order to check
the constants. [B

\begin{lemma} \label{del}An operator $T \in \L(X,Y)$ is of weak cotype
$g$ if and only if there is
a $0<\delta <1$ and a constant $C_{\delta}>0$ such that for all $\nen$
and $u \in \L(\lzn,X)$ one has
\[ g(n) \pl a_{[\delta n]}(Tu) \kl C_{\delta}(T) \pl \ell(u) \pl .\]
Moreover, we have the following relation for the constants
\[ \frac{\delta}{2S_2}  \pl C_{\delta}(T) \kl
wc_g(T) \kl e^{\frac{3}{2}}\pl S_2 \pl \kla 1-\delta
\mer^{-\frac{1}{2}} \pl C_{\delta}(T) \pl .\]
\end{lemma}

{\bf Proof:} The first estimate of $C_{\delta}$ by $wc_q(T)$ is
obvious.
For the second let $u \in \L(\ell_2,X)$, $v\in \Pi_2(Y,\ell_2)$ and
$\nen$.
By Schmidt decomposition there is subspace $H\subset \ell_2$ with $dim
H\lell n$
such that $a_n(vTu) \lell a_n(vTu\iota_H)$. We set $m\p:=\p n-[\delta
n]$ provided
$\delta n \gll 1$ and $m\p:=\p n$ else. Using the multiplicativity of
the
Weyl numbers and the Weyl number estimate
for the 2-summing norm we obtain
\for
a_n(vTu) \lel a_n(vTu\iota_H) &\le& a_{[\delta n]}(Tu\iota_H) \pl
x_m(v)
\kl C_{\delta} \pl g(n)^{-1}\pl \ell(u) \pl m^{-\frac{1}{2}} \pl
\pi_2(v) \\
&\le& C_{\delta}\pl (1-\delta)^{-\frac{1}{2}}\pl
(g(n)n^{\frac{1}{2}})^{-1}\pl \ell(u)\pl \pi_2(v)
\mel
From \cite[...]{DJ1} we deduce
\for
 wc_g(T) \kl S_2 \pl e^{\frac{3}{2}} \pl (1-\delta)^{-\frac{1}{2}} \pl
 C_{\delta} \pl .\\[-1.3cm]
\mel\qed

Now we can prove the proposition which initialize Talagrand's
machinery.

\begin{prop} \label{ini}Let $D\p:=\p
2^{\frac{9}{2}}e^{\frac{3}{2}}S_2^2$ and
$T \in \L(X,Y)$. If there
is a constant $c>0$ such that for all
$n \in \nz$, all vectors $x_1,..,x_n$ the condition
\[  \summ_1^n \bet \langle y^*,Tx_j \rangle \rag^2 \kl \noo y^* \rrm
\quad \mbox{and} \quad \noo Tx_j \rrm \gl \frac{1}{D} \quad \mbox{for
all} \quad j=1,..,n \]
implies
\[ \kla \ez \noo \summ_1^n g_j \p x_j \rrm^2 \mer^{\frac{1}{2}}
\gl \frac{1}{c} \pl g(n)  \pl, \]
then $T$ is of weBak cotype $g$ with
\[ wc_g(T) \kl c \pl .\]
\end{prop}

{\bf Proof:} There is no loss of generality to assume $T$ of finite
rank.
Indeed, if we can prove the assertion for all $T_{|E}$, $E$ a finite
dimensional
subspace of $X$ we obtain
\[ wc_q(T) \kl \sup_{E f.d.} wc_q(T_{|E}) \kl c\pl . \]
If $T$ is of finite rank we deduce form the lemma \ref{del} above that
there is a
positive real number $A$ with
\[ \frac{wc_q(T)}{2a} \pl < \pl A \pl < \pl C_{\frac{1}{2}}(T) \pl ,\]
where $a \p:=\p \sqrt{2}e^{\frac{3}{2}}S_2$. By definition there is an
$m \in \nz$ and an
operator $u \in \L(\ell_2^m,X)$ such that
\[ a_{[\frac{m}{2}]}(Tu) \pl > \frac{1}{16 S_2 a}
\quad \mbox{and} \quad \ell(u) \kl \frac{ g(m)}{ 16 S_2 aA } \pl. \]
We define $n \p:=\p [\frac{m}{4}]$. By definition of the weak cotype
$g$
we have
\[ g(n) \pl a_n(Tu) \kl wc_g(T) \pl \ell(u)
\kl \frac{ wc_q(T) \p g(m)}{16 S_2 aA } \pl <\pl \frac{1}{8S_2} \p
g(m)\pl. \]
Let us first assume $m \gll 4$, hence $n \gll \frac{m}{8}$. Since the
approximation numbers coincide with the Gelfand
numbers for operators on Hilbert spaces, there exists a subspace $H \su
\ell_2^m$
with $codim(H)<n$ such that
\[ \noo Tu_{|H} \rrm \kl \frac{ \frac{g(m)}{8S_2} }{g(n) a_n(Tu)}
\pl a_n(Tu) \kl \frac{1}{8S_2}\pl \frac{g(m)}{g(n)}
\kl 1 \pl .\]
From elementary pr[Boperties of the approximation
numbers we deduce with $codim(H)<n$
\[ \frac{1}{16S_2a} \pl < \pl  a_{[\frac{m}{2}]}(Tu) \kl a_{2n}(Tu) \kl
a_n(Tu_{|H}) \pl .\]
By a lemma probably due to Lewis, \cite{PIE}, there is an orthonormal
sequence
$(w_j)_1^n \in H$ such that
\[ \noo Tu(w_j) \rrm \gl \frac{1}{16 S_2a} \quad \mbox{for all} \quad
j=1,..,n \pl .\]
We define $x_j \p :=\p u(w_j)$ and $z \p :=\p \summ_1^n e_j \otimes x_j
: \ell_2^n \nach X$.
Since the system $w_j$ is an orthonormal sequence we get
$\noo Tz \rrm \kll \noo Tu_{|H} \rrm \kll 1$.
By the properties of the $\ell$ norm we deduce
\[ \kla \ez \noo \summ_1^n g_j x_j \rrm^2  \mer^{\frac{1}{2}}
\lel \ell(z) \kl \ell(u) \kl \frac{g(m)}{16S_2aA} \kl
\frac{g(n)}{2aA}   \pl .\]
By assumption this implies $2aA \kll c$ and therefore
\[ wc_g(T) \kl 2aA  \kl  c \pl .\]
If $m\kll 4$ we see that
\[ \noo Tu \rrm > \frac{1}{16S_2a}\quad \mbox{and}\quad \ell(u) \kl
\frac{1}{4aA}\pl .\]
From the weak cotype $g$ definition we have
\[ \noo Tu \rrm \kl wc_g(T) \pl \ell(u) <  \frac{g(m)}{(S_2} \kl
\frac{1}{2} \pl .\]
Let $h$ be a norm 1 vector were the norm is attained and let $x_1
\p:=\p 2\p u(h)$.
From $\noo Tx_1 \rrm \kll 1$, $\noo Tx_1 \rrm \gll \frac{1}{8S_2a}$ and
\[ \kla \ez \noo g_1 x_1 \rrm^2 \mer^{\frac{1}{2}}  \lel 2\pl
\ell(u_{|span\{h\}})
\kl \frac{1}{2aA}\]
we infer $2aA \kll c$ which implies the assertion.\qed

A precise calculation in the case $g(k) \lell k^{\frac{1}{q}}$, $2\kll
q <\infty$
shows that we can take $D\lell 8e$. Although the condition in the
proposition above is a little
bit technical it is nonetheless equivalent to the usual definition
of weak cotype $g$.[B

\begin{prop}\label{com}
If  $\p T\in \L(X,Y)$ is of weak cotype $g$,
$0<\rho <1$ and  $x_1,..,x_n \in X $ satisfying
\[ \sup_{y^* \in B_{Y^*}} \summ_1^n \bet \langle y^*, Tx_j \rangle
\rag^2 \kl 1
\quad \mbox{and} \quad \noo Tx_j \rrm \gl \rho \quad \mbox{for} \quad
j=1,..,n\]
then one has
\[ \rho^4 \p g(n) \kl S_2 \pl 2048 \pl wc_g(T) \pl
\kla \ez  \noo \summ_1^k g_j \p x_j \rrm^2 \mer^{\frac{1}{2}} \pl .\]
\end{prop}

{\bf Proof:} We choose functionals $(y_j^*)_1^n \subset B_{Y^*}$ with
\[ \rho \kl \noo Tx_j \rrm \lel \langle y_j^*,Tx_j \rangle \]
for $j=1,..,n$. We define the operators
\[ u\pl :=\pl \summ_1^n e_j \otimes x_j : \ell_2^n \nach X \quad
\mbox{and} \quad v\pl :=\pl \summ_1^n y_j^* \otimes e_j : Y \nach
\ell_{\infty}^n \pl .\]
By definition of $\noo v \rrm \kll 1$ and by assumption
$\noo Tu \rrm\kll 1$. On the other hand fix $1\kll k\kll n$.
Using an a well known estimate of the $2$-summing norm
by approximation numbers, \cite{PIE} we get
\for
\rho\pl \sqrt{n} &\le& \pi_2(vTu) \kl 2 \pl \summ_1^n
\frac{a_j(vTu)}{\sqrt{j}}\\
&\le& 2 \summ_1^{k-1} \frac{1}{\sqrt{j}} \pl +\pl 2\p
\frac{n}{\sqrt{k}} \pl a_k(vTu)
\kl 4 \sqrt{k-1}\pl + \pl  2\p \frac{n}{\sqrt{k}} \pl a_k(Tu) \pl.
\mel
If we choose $k-1 \kll \frac{\rho^2}{64} n \kll k$ we deduce
\[\rho \lel \frac{2}{\sqrt{n}} \pll \frac{\rho \sqrt{n}}{2}
\kl \frac{2}{\sqrt{n}[B} \pll  2 \p \frac{n}{\sqrt{k}} \pl a_k(Tu)
\kl \frac{32}{\rho} \pl a_k(Tu) \pl .\]
Finally we get
\for
g(n) &\le& \frac{S_2 64}{\rho^2} \pl g(k)
\kl S_2 \p \frac{2048}{\rho^4} \pl g(k) \pl \frac{\rho^2}{32}
\kl  S_2 \p \frac{2048}{\rho^4} \pl g(k) \pl a_k(Tu)
\kl  S_2 \p  \frac{2048}{\rho^4} \pl wc_g(T) \pl \ell(u)\pl .
\\[-1.5cm]
\mel \qed

Now we will briefly proof the easy implications of our theorem.

{\boldmath $ii) \Rightarrow iii)$\unboldmath} \pl
Let $T\lell RS$, where $R:X \nach \ell_{\infty}$ and $S:\ell_{\infty}
\nach X$.
By proposition \ref{com} $X$ is of weak cotype $g$ and in particular of
finite cotype
by condition $L)$.
From Maurey's theorem,
see \cite{TOJ} every operator $S: \ell_{\infty} \nach X$ is $p$-summing
for some $p <\infty$.
Using the gaussian version of Kintchine's inequality we get for all $u:
\ell_2 \nach \ell_{\infty}$
\for
\sup_k g(k) \pl a_k(Su) &\le& wc_g(X) \pl \ell(Su)
\kl wc_g(X)\pl \pl \sqrt{p}\pl \pi_p(S) \pl \noo u \rrm \\
&\le& wc_g(X)  \pl \sqrt{p} \pl c(p,X)  \pl \noo S \rrm \pl \noo u \rrm
\pl .
\mel
This means $\lgx(S) \kll \sqrt{7} \p c(p,X) \p wc_g(X) \p\noo S\rrm$
for some constant
$c(p,X)$. (Actually this part of the proof is due to Mascioni.)
To conclude we only have to note that $\Lgx$ is of eigenvalue type
$\ell_{g,\infty}$ by the generalized
Weyl's inequality, see \cite{PII}, together with condition $S_4$. Hence
we get
\for
\sup_k g(k)\pl \bet \lambda_k(T) \rag
&\le& c_1 \pl \lgx(T) \kl c_1 \noo R \rrm \pl \lgx(S)
\kl c_1 \pl \sqrt{p}\pl c(p,X) \pl wc_g(X) \pl \noo R \rrm \pl \noo S
\rrm \pl .\\[-1.3cm]
\mel\qed

{\boldmath $iii) \Rightarrow i)$\unboldmath}
Let $E$ be a $n$-dimensional subspace of $X$ and $(x_k)_1^N\subset E$.
We choose functionals $x_1^*,.., x_N^*$ in $B_{X^*}$
such that
\[ \noo x_i \rrm \lell \langle x_i, x_i^*\rangle \pl .\]
We define the operators
\[ v\pl := \pl \summ_1^N x_i^* \otimes e_i : X \nach \ell_{\infty}^N
\quad \mbox{and} \quad u\pl :=\pl \summ_1^N e_i \otimes x_i :
\ell_{\infty}^N \nach X \pl .\]
By definition of $v$ we have $\noo v \rrm\kll 1$ and
\[ \noo u \rrm \lell \sup_{\bet \alpha_i \rag \le 1} \noo \summ_1^N
\alpha_i x_i \rrm\pl .\]
Therefore we conclude
\for
\summ_1^N \noo Tx_i \rrm &=& tr(uv) \kl \summ_1^n \bet \lambda_k(uv)
\rag \\
&\le& \summ_1^n g(k)^{-1}  \pl \sup_k g(k) \pl \bet \lambda_k(uv)
\rag\\
&\le& S_4\pl \frac{n}{g(n)} \pl c_3 \pl \noo v \rrm \pl \noo u \rrm \\
&\le& 4\pl S_4 \pl c_3 \pl \frac{n}{g(n)} \pl \sup_{\alpha_i\pm 1}\noo
\summ_1^n \alpha_i x_i \rrm \\[-1.3cm]
\mel\qed

The implication {\boldmath $i) \Rightarrow {\it i'})$\unboldmath}
follows
obviously from the contraction principle, see \cite{LTII}
\[ \sup_{\bet \alpha_k \rag \le 1 }\noo \summ_1^n \alpha_k x_k \rrm
\kl 4 \sup_{\eps_k=\pm 1 }\noo \summ_1^n \eps_k x_k \rrm\]
and the trivial observation that $n$ elements are contained in a $n$
dimensional
subspace of $X$.\qed

Till the end of this chapter we are concerned with the proof of the
implication
\boldmath $ {\it i'}) \Rightarrow ii)$ \unboldmath. Assuming ${\it
i'})$ we first
observe
\for
\noo (\noo Tx_k \rrm)_k \rrm_{g,\infty} \kl c_1 \sup_{x^* \in B_{X^*}}
\summ_k \bet \langle x^*,x_k \rangle \rag \pl .\\[-1.5cm]
\mel \hfill $\bf (*)$\hz

Indeed this is a classical argument. We can assume $\noo Tx_i \rrm$ non
increasing and fix $k\in \nz$. Then we have
\for
k \noo Tx_k \rrm &\le& \summ_1^k \noo Tx_i \rrm \kl c'_1\pl
\frac{k}{g(k)} \pl  \sup_{\eps_k=\pm 1}\noo \summ_1^k \eps_i x_i \rrm
\\
&\le& c'_1 \pl \frac{k}{g(k)} \pl \sup_{x^* \in B_{X^*}}\p \sup_{\bet
\alpha_i \rag \le 1} \bet \langle x^*,\summ_1^k \alpha_i x_i \rangle
\rag\\
&\le& c'_1 \pl \frac{k}{g(k)} \pl \sup_{x^* \in B_{X^*}} \summ_k \bet
\langle x^*,x_i \rangle \rag \pl .
\mel
The best possible constant in $(*)$ will  be denoted by $H$. Let us
note that
for vectors $(x_i)_1^n$ with $\noo x_i \rrm \gll 1$ we certainly have
\[g(n)\kl H\pl \sup_{x^* \in B_{X^*}} \summ_1^n \bet \langle
x^*,x_i\rangle \rag \pl .\]
With this observation the following two lemmata from Talagrand can be

formulated in our setting. The first one is a lemma which allows to
regroup a certain collection of disjoint blocs.

\begin{lemma}[Talagrand: Lemma 4.2.]\label{regroup} There exists a
constant $K>0$ with the following property.
Consider disjoint subsets $I_1,..,I_k$ of $\{1,..,n\}$ with union $I$.
Let $\alpha > 0$ such that for all
$1\kll j \kll k$
\[ \ez \noo \summ_{i \in I_j} g_i x_i \rrm \gl \alpha
\quad \mbox{and} \quad
\sqrt{k} \kl \frac{\alpha }{2S_3KH}\pl g(k) \pl  .\]
Then one has
\[ \ez \noo \summ_{i \in I} g_i x_i \rrm \gl \frac{\alpha}{2S_3H} \pl
g(k) \pl .\]
\end{lemma}

In the following we want to prove that $(*)$ implies weak cotype $g$.
By proposition \ref{ini}
and Kahane's inequality we are left to verify that for all vectors
$x_1,..,x_n$ with $\noo x_i \rrm \gl \frac{1}{D}$
and $\summ_1^n \bet \langle x^*,x_i \rangle \rag^2 \kll \noo x^* \rrm$
one has
\[ \ez \noo \summ_1^n g_i x_i \rrm \gl \frac{g(n)}{c(g,H)} \pl. \]
Therefore we will fix in the following this sequence of vectors. The
second
lemma proved with the concentration of measure phenomena reads as
follows
[B
\begin{lemma}[Talagrand: Lemma 4.3.] \label{select}Let $8\kll s\kll n$
such that
\[ \frac{s}{\sqrt{n}} \kl \frac{1}{16HD} \pl g(s)\]
and $J$ a subset of $\{1,..,n\}$ with $cardJ\gll \frac{n}{2}$.
Then there exists a subset $I\su J$ with $card I\lell s$ and
\[ \ez \noo \summ_{i \in I} \eps_i x_i \rrm \gl \frac{g(s)}{64HD}\pl
.\]
\end{lemma}

\begin{rem} The most interesting application of the theorem
is certainly given by the sequence $g(k) \lell \sqrt{k}$. In this case
the conclusion of the theorem is very easy. Indeed, we can choose
$s\kll \frac{n}{ 64H^2D^2} \kll 2s$ and get a sequence $I \subset
\{1,..,n\}$
with $cardI\lell s$ such that
\[ \ez \noo \summ_1^n g_i x_i \rrm
\gl \sqrt{\frac{2}{\pi}} \pl \ez \noo \summ_{i \in I}\eps_i x_i  \rrm
\gl \sqrt{\frac{2}{\pi}}\pl \frac{\sqrt{s}}{64HD}
\gl \sqrt{\frac{2}{\pi}} \pl \frac{1}{(64HD)^2} \pl \sqrt{n} \pl .\]
Proposition \ref{ini} implies
\[ wc_2(X) \kl c_0 \pl H^2 \pl .\]
\end{rem}

Now we start with the main proof. In the sequel we will assume some
conditions
to be verified. At the end we will discuss the influence of this
conditions
on the constants. We choose an even number $M \in \nz$ such that
\[ 2^{rM+1} \kl n \kl 2^{rM+1} 2^{2r}\pl .\]
Si[Bnce $ g(2^{2r}2^{rM+1}) \kll S_2 \p 2^{2r} \p g(2^{rM+1})$
we can even assume $n\lell 2^{rM+1}$.
Furthermore, we set
\[ N \pl:=\pl \frac{M}{2}\p ,\quad  p \pl :=\pl s \pl :=\pl 2^{rN} \pl
.\]
The condition
\for
\sqrt{2}\p 16\p HD \kl g(s)\quad \mbox{and}\quad 8 \kl s  \\[-1.5cm]
\mel \hfill $\bf (1)$
\hz

implies $\frac{s}{\sqrt{n}} \kll \frac{g(s)}{16HD}$. From a successive
application of lemma \ref{select} we can find $p$ disjoint subsets
$I_1,..,I_p$ each of cardinality $s$ such that
\for
\ez \noo \summ_{i \in I_j} g_i x_i \rrm &\ge& \sqrt{\frac{2}{\pi}} \pl
\ez \noo \summ_{i\in I_j} \eps_i x_i \rrm \gl  \frac{g(s)}{100 HD} \pl
.\\[-1.5cm]
\mel\hfill $\bf (\circ)$\hz

Now we will apply the iteration procedure to regroup disjoint blocs.
This lemma is also essentially contained in \cite{TAL}.

\begin{lemma} \label{iterate} Let $k  \in \nz$ satisfy
\for
g(k)\gl 2S_3(K+1)H \quad \mbox{and}\quad g(s)\gl 100\p D  \p H  \p
\sqrt{k}\pl .  \\[-1.5cm]
\mel \hfill $\bf (2)$
\hz

Given a subset $T \subset \{1,..,p\}$ with $card(T) \lell k^l$,
$k^l\kll p$
and $I_T \lell \cup_{j\in T} J_t$ one has
\[ \ez \noo \summ_{i \in I_T} g_i x_i \rrm \gl \frac{g(s)}{100DH} \pl
\kla \frac{g(k)}{2S_3H} \mer^l \pl . \]

\end{lemma}

{\bf Proof:} The case $l=0$ is $(\circ)$. Proceeding
by induction we can assume that the statement is valid for $l$. A set
$T$ of cardinality $k^{l+1}$ can be split up into $k$ sets $T_j$ with
cardinality
$k^l$. By induction hypothesis we have for all $j =1,..,k$
\[ \ez \noo \summ_{i \in I_{T_j}} g_i x_i \rrm \gl \frac{g(s)}{100DH}
\pl \kla \frac{g(k)}{2S_3H} \mer^l  \pl =: \pl \alpha \pl .\]
The assertion follows from lemma \ref{regroup} provided we have
\[ \sqrt{k} \kl \frac{\alpha}{2S_3KH} \pl g(k)
\lel \frac{g(s)}{100HD}\pl \frac{g(k)}{2S_3KH} \pl \kla
\frac{g(k)}{2S_3H_3} \mer^l   \pl .\]
which is obvious by our assumption. \qed

Now we set $k\p:=\p 2^N$ and assume $(1)$ and $(2)$ to be satisfied.
Then we have
can apply lemma \ref{iterate} to find
\[ \ez \noo \summ_{i \in \cup_{j=1,..,p} J_j} x_i g_i \rrm
\gl \frac{1}{100DH}\pl \kla \frac{1}{2S_3H} \mer^r \pl g(s) \p g(k)^r
\gl \frac{1}{100DM_rH}\pl \kla \frac{1}{2S_3H} \mer^r \pl g(2^{2Nr})
\pl .\]
In this case we set
\[ c_1(g,H) \pl:=\pl S_2 \pl 2^{2r+1} \pl 100 DM_r \pl (2S_3)^r \pl
H^{r+1} \pl .\]

In order to garuantee the conditions $(1)$ and $(2)$ we define
\[ B\pl :=\pl \max\{2S_3L_t(K+1), 100L_tD\}^{\max\{2,t\}} \]
and assume first $k \gll BH^{\max\{2,t\}}$. Since $H\gll 1$, $L_t\gll
1$,
we trivially have $s\gll 8$. Furthermore, we get
\[ 2S_3(K+1)H \kl \frac{k^{\frac{1}{t}}}{L_t} \kl g(k)
\quad\mbox{and}\quad
100DH\p \sqrt{k}\kl \frac{\sqrt{k}}{L_t} \p \sqrt{k}
\kl g(k^{r}) \lell g(s) \pl .\]
If the remaing case $k \kll BH^{\max\{2,t\}}$ we define
\[ c_2(g,H) \pl :=\pl S_2 \pl 2^{2r+3}\pl B^r \pl
g([H^{\max\{2r,tr\}}]) \pl .\]
For $c(g,H) \lell \max\{ c_1(g,H),c_2(g,H)\}$ we have $\ez \noo
\summ_1^n g_i x_i \rrm \gll \frac{g(n)}{c(g,H)}$
in any case and the proof is finished.
\begin{rem}
\begin{enumerate}
\item In the case $g(k) \lell k^{\frac{1}{q}}$ the prove above
gives a constant of order $c_q \p H^{2q+2}$ which is certainly not
optimal.
\item If the condition $M)$ is not satisfied we can define the new
sequence
\[ \tilde{g}(n) \pl :=\pl \max\{ g(k^r)g(k)^r \p \mitt \p k^{2r} \kll n
\} \pl \]
It is easy to check that the conditions $Si)-Sii)$ as well as $L)$
are still satisfied (for probably different constants).
The proof above shows that a summing condition of order $\lgi$ for a
Banach space $X$
implies weak cotype $\tilde{g}$.
\end{enumerate}
\end{rem}

\section{Optimal summing and cotype spaces}
\setcounter{lemma}{0}

In the following we will define sequence spaces which are associated
with the cotype and summing properties of Banach spaces. In this
setting
it is more convenient to study the Rademacher Cotype. We
will use the following definition of a maximal symmetric sequence
space $Y$ which is a sequence space with the following properties

\begin{enumerate}
\item[i)] $\noo \tau \rrm_{\infty} \kll \noo \tau \rrm_Y \kll \noo \tau
\rrm_1$ for all sequences with finite support.
\item[ii)] $\noo \tau^*\rrm \lel \noo \tau \rrm$\p, where $\tau^*$
denotes
the non increasing rearrangement of $\bet \tau \rag$.
\item[iii)] $\noo \tau \rrm \lel \sup_n \noo P_n(\tau) \rrm$\p,
where $P_n$ denotes the projection onto the first $n$ coordinates.
\end{enumerate}

An operator $T\in \L(X,Y)$ is said to be $(Y,1)$-summing, of cotype
$Y$, if there is
a constant $c>0$ such that for all $\nen$, $x_1,..,x_n$ one has
\[ \noo \summ_1^n \noo  Tx_k \rrm \p e_k \rrm_Y \kl c \pl \sup_{x^*\in
B_{X^*}} \summ_1^n \bet \langle x^*,x_k \rangle \rag \p, \quad
 \noo \summ_1^n \noo  Tx_k \rrm \p e_k \rrm_Y \kl c \pl \ez \noo
 \summ_1^n \eps_k \p x_k \rrm \p ,\pl {\rm resp.}\]
The corresponding norm is denoted by $\pi_{Y,1}(T)\p:=\p \inf\{c\}$,
$c_Y(T) \p:=\p \inf\{c\}$,
where the infimum is taken over all $c$ satisfying the inequality
above.

\begin{rem} \label{ws}\rm This definition can in particular
be applied for $Y \lell \ell_{g,infty}$. We want to compare cotype
$\ell_{g,\infty}$ with the notion of weak (gaussian) cotype $g$ defined
in the chapter before. Using Lewis' Lemma
as in the proof of proposition \ref{ini} it is quite easy to see that
Rademacher cotype $\ell_{g\infty}$ implies weak cotype $g$.
The converse is not always true. If we consider $g(n)\lell \sqrt{n}$
we see that every Banach spaces with weak cotype 2, not having cotype
2,
yields an example of a Banach space having cotype $g$, but not
cotype  $\ell_{g\infty}$. A further example is given by
$g(n) \lell \sqrt{1+\ln n}$, since every Banach space has
cotype $g$ but cotype $\ell_{g\infty}$ only holds for Banach spaces
with
finite cotype by Maurey/Pisier's theorem. Therefore,  it is natural
to require two additional conditions, namely $g(n)\gl c\p
n^{\frac{1}{p}}$
for some [B$p < \infty$ and
\[ g(n) \kl c_q \pl \kla \frac{n}{k}\mer^{\frac{1}{q}} \pl g(k) \]
for some $q>2$. In this case a Banach space with cotype $g$ is of
finite
cotype and using the inequality
\for
\pi_2(T) &\le& C_q\pl c_q \pl \frac{ \sqrt{n} }{ g(n) } \pl \sup_k
g(k)\pl x_k(T) \\
\mel
valid for all operators of rank at most $n$, we easily see
that every Banach space of cotype $g$ is also of cotype
$\ell_{g,\infty}$.
\end{rem}

The main tool of this chapter are the properties
of the optimal summing and cotype space associated to a Banach space
$X$.
Given  $\tau \lell (\tau_k)_k$ we define
\[ \tno \alpha \tno_S \pl :=\pl \inf\left \{ \sup_{| \alpha_k| \le
1}\noo \summ_1^n \tau_k \alpha_k x_k \rrm \mitt (x_k)_1^n \subset X,
\pl \noo x_k \rrm \lell 1 \right \} \]
and
\[ \tno \alpha \tno_{C} \pl :=\pl \inf\left \{ \ez \noo \summ_1^n
\eps_k \alpha_k x_k \rrm \mitt (x_k)_1^n \subset X, \pl \noo x_k \rrm
\lell 1 \right \} \pl .\]
Clearly this are homogeneous expression which are invariant under
permutations
and change of signs. In order to guarantee the triangle inequality we
define for
$T\in\{C,S\}$
\[ \noo \tau \rrm_T^0 \pl :=\pl \inf \left\{ \summ_1^m \tno \tau^j
\tno_T \mitt  n\in \nz, \p \tau^j \pl \mbox{with finite support
and}\quad \bet \tau \rag \kll \summ_1^m \bet \tau^j\rag \right \}\]
and
\[\noo \tau \rrm_T \pl :=\pl \sup_n \noo P_n(\tau) \rrm_{T}^0 \pl .\]
The space $Y_S$, $Y_C$ defined by this norm
will be called {\it optimal summing space, optimal cotype space},
respectively.
We summarize the properties of this spaces in the following

\begin{lemma} \label{opt}Let $X$ be a Banach space and $Y_S$, $Y_C$
it's optimal summing, optimal
cotype space, respectively, and let $Z$ be a maximal sequence space.
Then one has
\begin{enumerate}
\item The identity of $X$ is $(Y_S,1)$ summing and of cotype $Y_C$ with
constant $1$.
\item The identity of $X$ is $(Z,1)$-summing (of cotype $Z$) if and
only if
\[ Y_S \pl \subset \pl Z\quad (Y_C \pl  \subset \pl Z, \quad resp.)\]
The norm of the inclusion is $\pi_{(Z,1)}(id_X)$ ($C_Z(id_X)$,
respectively).
\item For $Y\in \{Y_S,Y_C\}$ and each finitely supported sequence
$(\tau^k)_1^n$
one has
\[ \noo \summ_1^n \noo \tau^k \rrm_Y \p e_k  \rrm_Y \kl \noo \summ_1^n
\bet \tau^k \rag \rrm_Y \pl .\]
\end{enumerate}
\end{lemma}

{\bf Proof:} $1.$, $2.$ are obvious. We will only consider
the cotype case in $3.$.
We denote by $\tau \p:=\p \summ_k \bet \tau^k \rag$. Given $\delta >0$
we can find a finite sequence $(x_i)_i\subset X$ with $\noo x_i \rrm
\lell 1$
such that
\[ \ez \noo \summ_i \eps_i \tau_i x_i \rrm \kl (1+\delta)\pl \tno \tau
\tno_C \pl .\]
For any sequence of signs $\rho_k$ we can find a sequence
$(\gamma_i)_i$, $\gamma_i \in [-1,1]$
such that
\[ \summ_k \rho_k \bet \tau^k \rag \lell \gamma \tau \pl .\]
By the sign invariance of $(\eps_i)$ and the fact that extreme points
in the unit ball $\lin$ over $\rz$ are sequences of signs, see e.g.
\cite{PIE},
we get
\[  \ez_{\eps} \noo \summ_i \eps_i \kla \summ_1^n \rho_k \bet \tau^k_i
\rag \mer x_i \rrm \lel
 \ez_{\eps} \noo \summ_i \eps_i \gamma_i \tau_i x_i \rrm \kl
 \ez_{\eps} \noo \summ_i \eps_i \tau_i x_i \rrm \kl (1+\delta) \pl \tno
 \tau \tno_C \pl .\]
Taking expectations we deduce from $1.$ and the triangle inequality in
$Y$
\for
(1+\delta) \pl \tno \tau \tno_C &\ge&
 \ez_{\eps}B \ez_{\rho} \noo \summ_i\eps_i \kla \summ_k \rho_k \bet
 \tau^k_i\rag \mer x_i \rrm
\gl\ez_{\eps} \noo \summ_k \noo \summ_i \eps_i \bet \tau_i^k\rag x_i
\rrm_X \p e_k \rrm_{Y_C}\\
&\ge& \noo \summ_k \kla \ez_{\eps} \noo \summ_i \eps_i \bet \tau^k_i
\rag  x_i \rrm \mer\p e_k \rrm_{Y_C}
\gl \noo \summ_k \tno \tau^k\tno_C \p e_k \rrm_{Y_C}\pl .
\mel
Letting $\delta$ to zero we have proved
\for
\noo \summ_k \tno \tau^k \tno_C \pl e_k \rrm_{Y_C} \kl \tno \summ_k
\bet \tau^k \rag \tno_C \pl .\\[-1.5cm]
\mel \hfill $\bf (*)$\hz

Now let $\tau \kl \summ_j \bet \si^j\rag$. We define $\beta^k \p :=\p
\frac{1}{\bet\tau^k\rag}\tau$
by pointwise multiplication and using the convention $\frac{0}{0}\lell
0$.
For the sequences $\si^{kj} \p:=\p \beta^k \si^j$ we clearly have
\[ \summ_k \bet \si^{kj}\rag \kl \bet \si^j \rag \quad \mbox{and}
\quad \bet \tau^k \rag \kl \summ_j \bet \si^{kj} \rag \pl .\]
From $(*)$ apllied for each sequence $(\bet \si^{kj}\rag)_k$ we deduce
\for
\noo \summ_k \noo \tau^k \rrm_Y \p e_k \rrm_{Y_C} &\le&
\noo \summ_k \kla \summ_j \tno \si^{kj} \tno_C \mer \p e_k\rrm_{Y_C}
\kl  \summ_j \noo \summ_k \tno \si^{kj}\tno \p e_k \rrm_{Y_C}\\
&\le& \summ_j \tno \summ_k \bet \si^{kj}\rag \tno_C \kl
\summ_j \tno \si^j \tno_C\pl .
\mel
Taking the infimum over all $\tau \kll \summ_j \si^j$ yields the
assertion.
\qed

We will study in more detail the spaces which satisfy the last
condition.
For a maximal symmetric sequence space $Y$ we define
\[ f_Y(n) \pl :=\pl \noo \summ_1^n e_k \rrm_Y \quad \mbox{and}\quad
q_Y \pl:=\pl \inf\left \{ 0<q<\infty \mitt \exists\p C: \pll
n^{\frac{1}{q}} \kl C \pl f_Y(n)  \right \}\pl .\]
Obviously we have the following inclusions for all $q>q_Y$
\[ Y\pl \subset \pl \ell_{f_Y,\infty} \pl \subset \pl \ell_q \pl .\]
If $Y$ satisfies the concavity condition $3.$  we have the following
alternative which is somehow an improvement of the classical
Maurey/Pisier argument in the context of finite cotype.

\begin{prop} \label{tensor}Let $Y$ be a maximal symmetric sequence
space which satisfies
\[ \noo \summ_1^n \noo \tau^k \rrm_Y \p e_k  \rrm_Y \kl \noo \summ_1^n
\bet \tau^k \rag \rrm_Y \pl .\]
For all $1 \kll p < \infty$ we have either $\ell_p \subset Y$ with
inclusion norm $1$ or there
exists a $q<p$ with $Y\subset \ell_q\p {\subset \atop \neq}\p \ell_p$.
In particular, we have
\[ \ell_{q_Y} \pl \subset \pl Y \pl .\]
\end{prop}

{\bf Proof:} Let $\tau$ be a sequence of finite support, $\tau_k \lell
0$ for
$k \gll n$ say. For $i=1,..,n$ we define
\[ \si_i \pl :=\pl \summ_{j=1}^{n} \pl \tau_j \pl e_{(i-1)n+j} \quad
\mbox{and the product}\quad
 \tau \otimes \tau \pl :=\pl  \summ_{i=1}^n \tau_i \pl \si_i \pl .\]
Usually, $\tau \otimes  \tau $ is defined in $\ell_p(\nz^2)$ which is
isometric
isomorph to $\ell_p$ by a renumbering of $\nz^2$. That's what we try to
imitate
with the definition above. Clearly, we have $\noo \tau \otimes\tau
\rrm_p \lell \noo \tau \rrm_p^2$.
Our assumption on $Y$ implies
\for
 \noo \tau \rrm_Y^2 &=& \noo \summ_1^n \tau_i \p \noo \si_i \rrm_Y \p
 e_i \rrm_Y
\kl \noo \tau \otimes \tau \rrm_Y \pl .\\[-1.5cm]
\mel \hfill $\bf(*)$\hz

In particular, $f_Y$ is submultiplicative, i.e. $f_Y(n) \p f_Y(k) \kll
f_Y(nk) \pl$.
Now we consider the following alternative.
\begin{enumerate}
\item There exists an $n_0 \in \nz$ such that  $f_Y(n_0)\pl >\pl
n_0^{\frac{1}{p}}$.
\item For all $n \in \nz$ one has $f_Y(n) \kl n^{\frac{1}{p}}$.
\end{enumerate}
In the first case  we choose $q<p$ such that $f_Y(n_0) \lell
n_0^{\frac{1}{q}}$.
For $n \in \nz$ we choose $m\in\nz$ with $n_0^{m-1} \kll n \kll
n_0^m$.
By the submultiplicativity and the triangle inequality we deduce
\for
n^{\frac{1}{q}} &\le& n_0^{\frac{m}{q}} \lel  f_Y(n_0)^m
\kl f_Y(n_0^m) \kl n_0 \pl f_Y(n) \pl .
\mel
This means $Y \subset \ell_{q,\infty}$ and therefore for all $q<r<p$
the inclusion
$Y \p \subset\p  \ell_r \pl {\subset \atop \neq}\pl \ell_p$.
Now we consider the second case.
We will first show $\ell_{p,1} \subset Y$.
Indeed, let $\tau$ ba a non increasing positive sequence with finite
support.
Then we have
\for
\noo \tau \rrm_Y &\le& \summ_{k=0}^{\infty} \noo
\summ_{j=2^k}^{2^{k+1}} \tau_j e_j \rrm_Y
\kl \summ_{k=0}^{\infty} \tau_{2^k} f_Y(2^k)
\kl\summ_{k=0}^{\infty} \tau_{2^k} (2^k)^{\frac{1}{p}}\kl 5 \pl \noo
\tau \rrm_{p,1} \pl .
\mel
Defining   $C_n \p :=\p \noo P_n : \ell_p \nach Y \rrm$, this means
$C_n \kl 5(1+\ln n)$. Now we will use a tensor trick
to finish the proof. For this we prove $C_n^2 \kll C_{n^2}$
Indeed, let $\tau$ a sequence with support contained in $\{1,..n\}$.
From $(*)$ we deduce
\[ \noo \tau \rrm_Y^2 \kl \noo \tau \otimes \tau \rrm_{Y}
\kl C_{n^2} [B\pl \noo \tau \otimes \tau \rrm_p \kl C_{n^2} \noo \tau
\rrm_p^2 \pl .\]
Hence we get
\for
C_n &\le& \inf_k \kla C_{ n^{2^k} } \mer^{\frac{1}{2^k}}
\kl \inf_k \kla 5 (1+ 2^k \ln n) \mer^{\frac{1}{2^k}} \lell 1 \pl .\\
\mel
If we apply the alternative for the space $\ell_{q_Y}$ we only have
to observe that an inequality $Y\subset \ell_q$ implies
$n^{\frac{1}{q}} \kl C \pl f_Y(n)$. By definition
this is impossible for all $q<q_Y$.\qed

As an application we will investigate cotype properties with respect to
the Lorentz space $\ell_{q,w}$.

\begin{prop} Let $2\kll q<\infty$, $1 \kll w \kll\infty$. A Banach
space $X$ is of cotype
$\ell_{q,w}$ if and only if
\[ X \mbox{ is of cotype } \quad \left \{
\begin{array}{l@{\quad \mbox{if} \quad}l}
\mbox{$p$ for some }p<q & w<q\\
\ell_{p,\infty} & q<w \pl .\\
\end{array}\right.\]
If $X$ is of cotype $\ell_{q,\infty}$ there exists a constant $C$ such
that
\[ c_q(id_E) \kl \sqrt{\pi} \pl C^{k+1} \pl (\max\{1,{\rm
log}_2\})^{(k)}((1 +{\rm log}_2 n)^{\frac{1}{q}})\]
holds for all $k \in \nz$ and $n$-dimensional subspaces $E \subset X$.
In particular,
\[ c_q(id_E) \kl \sqrt{\pi}\pl 2 \pl  C^{1+k_n}\pl ,\]
where $k_n$ is the smallest integer $k$ with $n\kl \underbrace{
2^{2^{.^{.^2}}} }_{k {\rm times}}$.
\end{prop}

{\bf Proof:} If $w<q$ and $X$ is of cotype $\ell_{q,w}$ we have $Y_C
\subset \ell_{q,w}$
by lemma \ref{opt}, but certainly not $\ell_q \subset Y_C$. By lemma
\ref{tensor} there must be
a $p<q$ such that $Y_C\subset \ell_p$. Since $X$ is of cotype $Y_C$ it
is also of cotype $\ell_p$.

Now let us assume that $X$ is of cotype $\ell_{q,\infty}$ with constant
$D$, say.
This implies in particular $n^{\frac{1}{q}} \kl C\pl f_{Y_C}(n)$. Let
$q\kll w <\infty$
and $\tau$ a positive non increasing sequence of finite support. For
$k\in \nz$ we define the disjoint elements
$ x_k \p :=\p \tau_{2^k} \summ_{j=2^{k-1}+1}^{2^k} e_j$.
Then we get
\for
\noo \tau \rrm_{q,w} &\le& \kla \summ_n (\tau_n n^{\frac{1}{q}})^w \p
\frac{1}{n} \mer^{\frac{1}{w}}
\kl  \kla \summ_{k=0}^{\infty} (\tau_{2^k} \p 2^{\frac{k}{q}})^w
\mer^{\frac{1}{w}}
\kl \noo \tau \rrm_{\infty} + 2^{\frac{1}{q}}D \pl \noo \summ_{k\in
\nz} \noo x_k \rrm_{Y_C} \pl e_k \rrm_w \pl.
\mel
If $w>q$ we have $\ell_{q,\infty} \subset \ell_w$ with inclusion norm
$c_{qw}$.
Therefore we deduce from $\summ_k x_k \kl \tau$ and condition $3.$ in
lemma \ref{opt}
\for
\noo \tau \rrm_{q,w}&\le&
 2^{\frac{1}{q}}D \pl \kla \noo \tau \rrm_{\infty} \pl +\pl c_{qw}\pl
 \noo \summ_k \noo x_k \rrm_{Y_C} e_k \rrm_{q,\infty} \mer\\
&\le& 2^{\frac{1}{q}}D \pl \kla  \noo \tau \rrm_{Y_C} \pl +\pl c_{qw}\p
D\pl \noo \summ_k \noo x_k \rrm_{Y_C} e_k \rrm_{Y_C} \mer
\kl 2^{1+\frac{1}{q}}\p c_{qw}\p D^2 \noo \tau \rrm_{Y_C} \pl .
\mel
Since $X$ is of cotype $Y_C$ and $Y_C \subset \ell_{q,w}$ we obtain the
assertion in this case.
Now we come to the case $q=w$. We denote by $\alpha_n \p:=\p \noo P_n:
Y_C \nach \ell_q^n\rrm $, with
the convention $\alpha_0 \lell 1$. We will prove
\[ \alpha_n \kl 2^{1+\frac{1}{q}}D\pl \alpha_{[{\rm log}n]} \pl .\]
Indeed, if the support of the given sequence $\tau$ above is contained
in $\{1,..,n\}$
then we will have $x_k\lell 0$ whenever $2^k > n$. Therefore we obtain
again by lemma \ref{opt}
\for
\noo \tau \rrm_q
&\le& \kla \noo \tau \rrm_{\infty} + 2^{\frac{1}{q}}D \pl \noo
\summ_{k=1}^{[{\rm log}_2n]} \noo x_k \rrm_{Y_C} \pl e_k \rrm_q\mer \\
&\le& \kla \noo \tau \rrm_{\infty} + 2^{\frac{1}{q}}D \pl \alpha_{[{\rm
log}_2n]} \noo \summ_{k=1}^{[{\rm log}_2n]} \noo x_k \rrm_{Y_C} \pl e_k
\rrm_{Y_C} \mer
\kl (1+2^{\frac{1}{q}}D\alpha_{[{\rm log}_2n]}) \pl \noo \tau
\rrm_{Y_C} \pl .\\
\mel
Together with the trivial estimate $\noo id: \ell_{q,\infty}^n \nach
\ell_q^n \rrm \kll (1+{\rm log}_2 n)^{\frac{1}{q}}$
and induction this implies
\[ Rc_q^n(id_X) \kl D^{k+1} \pl (\max\{1,{\rm log}_2\})^{(k)}((1 +{\rm
log}_2 n)^{\frac{1}{q}})\pl .\]
Since the gaussian cotype constant of a $n$ dimensional space $E$ can
be well estimated
by the gaussian cotype constant, see \cite{TOJ,DJ2}, and as a
consequence of the comparism principle of
gaussian and Rademacher variables, see e.g. \cite{TOJ}, we deduce
\for
 c_q(id_E) &\le& \sqrt{2}\p c_q^n(Id_E) \kl \sqrt{\pi} Rc_q^n(Id_X)
\kl \sqrt{\pi} D^{k+1} \pl (\max\{1,{\rm log}_2\})^{(k)}((1 +{\rm
log}_2 n)^{\frac{1}{q}})\pl .\\[-1.5cm]
\mel \qed

\begin{frem} \rm
\begin{enumerate}
\item The same contraction argument can also be applied in the space
$Y_S$, provided we have $Y_S \subset \ell_{q,\infty}$. This is
only interesting in the case $q=2$. Hence in a weak cotype $2$ space we
have
\[ \pi_{21}^n(id_X) \kl C^{k+1} \pl (\max\{1,{\rm
log}_2\}^{(k)}((1+{\rm log}_2n)^{ \frac{1}{2} })\pl .\]
It is still open whether such an estimate is valid for the cotype $2$
constant.
\item Let $X$ be a Banach space with non-trivial cotype. Then $q_{Y_S}$
is finite
and we can fix a natural number $q_{Y_S} <r$. If we define
\[ g_S(n) \pl :=\pl f_{Y_S}([n^{\frac{1}{2}}]) \pl
f_{Y_S}([n^{\frac{1}{2r}}])^r \pl ,\]
we can apply the proof of $i) \Rightarrow ii)$ of the main theorem to
deduce
that $X$ is of weak cotype $g_S$. Given a number $q>2$ and a sequence
$g(n)_{n \in \nz}$ we define
\[ g_q(n) \pl :=\pl n^{\frac{1}{q}} \inf\{k^{-\frac{1}{q}}\p g(k) \mitt
1\kll k\kll n\} \pl .\]
From remark \ref{ws} we have clearly
\[ \frac{1}{c_{r,q}}\pl (g_S)_q(n) \kl f_{Y_C}(n) \kl f_{Y_S}(n) \pl
.\]
Nevertheless, the eigenvalue estimate for the sequence $(f_{Y_S})_q$
can be directly derived
from a generalization of Maurey theorem's theorem to $p$-convex
sequences, see \cite{J1,J2}.
(The space $\ell_{f_{q,Y_S},\infty}$ is $p$-convex for each $2<p<q$.)
This is of particular interest if the cotype properties of Orlicz
spaces associated to the function  $M(t) \lell \kla \frac{t}{1+\bet\ln
t\rag}\mer^{\frac{1}{q}}$
are studied in more detail.
\end{enumerate}
\end{frem}

\newcommand{\k}{\normalsize}

\hz
\k
1991 Mathematics Subject Classification: Primary:46B07, Secondary
47B06.
\begin{quote}
Marius Junge

Mathematisches Seminar der Universit\H{a}t Kiel

Ludewig-Meyn-Str. 4

24098 Kiel 1

Email nms006@rz.uni-kiel.d400.de

Germany
\end{quote}


\begin{thebibliography}{10001}
\bibitem[DJ1]{DJ1} M. Defant and M. Junge: {\k\sl On weak $r2$-summing
operators and weak Hilbert spaces;} Studia math. 96 (1990); 203-217.
\bibitem[DJ2]{DJ2} M. Defant and M.Junge: {\k\sl On absolutely summing
operators with apllication to the (p,q)-summing norm with few vectors};
J. of Functional Ana. 103 (1992), 62-73.
\bibitem[J1]{J1}  M. Junge: {\k\sl Comparing Rademacher and gaussian
cotype for operators on the space of continuous functions}; preprint.
\bibitem[J2]{J2}  M. Junge: {\k \sl Orlicz properties in operator
spaces and eigenvalue estimates}; preprint.
\bibitem[LTII]{LTII} J. Lindenstrauss and L. Tzafriri: {\k\sl Classical
Banach spaces II, function  spaces}; Springer Berlin Heidelberg New
York 1979.
\bibitem[MAS]{MAS} V. Mascioni: {\k \sl On weak cotype and weak type in
Banach spaces}; Note di Matematica Vol VIII-n.1(1988), 67-110.
\bibitem[MP]{MP} B. Maurey and G. Pisier: {\k\sl S\'{e}ries des
variables al\'{e}atoires vectorielles ind\'{e}ependentes et
 g\'{e}om\'{e}trie des espaces de Banach;} Studia Math. 58 (1976),
 45-90.
\bibitem[PI]{PIE} A. Pietsch: {\k \sl Operator Ideals}; Deutscher
Verlag Wiss.\p, Berlin 1978 and North Holland, Amsterdam-New
York-Oxford 1980. Cambridge University Press, 1987.
\bibitem[PII]{PII} A. Pietsch: {\k \sl Eigenvalues and s-numbers};
Cambridge university press, 1987.
\bibitem[PS]{PS} G. Pisier: {\k \sl Factorization of linear operators
and Geometry of Banach spaces;} CBMS Regional Conference Series
$n^{\circ}$ 60, AMS 1986.
\bibitem[PSW]{PSW} G. Pisier: {\k \sl Weak Hilbert spaces;} Proc.
London Math. Soc. 56(1988), 547-579.
\bibitem[TA]{TA} M. Talagrand : {\k\sl Orlicz property and cotype in
symmetric sequence spaces;} preprint.
\bibitem[TAL]{TAL} M. Talagrand: {\k\sl Cotype and $(q,1)$-summing norm
in a Banach space;} Invent. math. 110 (1992), 545-556.
\bibitem[TOJ]{TOJ} N. Tomczak-Jaegermann: {\k \sl Banach-Mazur
distances and finite-dimensional operator ideals}; Longmann, 1988.
\end{thebibliography}
\end{document}